\documentclass[12pt]{article}
\usepackage{mathtools}
\usepackage[utf8]{inputenc}
\usepackage{amssymb,amsmath,amsfonts,eucal,mathrsfs,amsthm} 
\usepackage[colorinlistoftodos]{todonotes}
\usepackage{xcolor}
\usepackage{enumitem}
\newtheorem{theorem}{Theorem}
\newtheorem{proposition}[theorem]{Proposition}
\newtheorem{lemma}[theorem]{Lemma}

\theoremstyle{definition}

\newcommand{\Z}{\mathbb{Z}}
\newcommand{\Q}{\mathbb{Q}}
\newcommand{\Sf}{\mathbb{S}}

\newcommand{\Hy}{\mathbb{H}}

\newcommand{\spa}{\mbox{span}}

\newcommand{\Ric}{\mbox{Ric}}

\newcommand{\trace}{\mbox{tr\,}}

\def\<{{\langle}}
\def\>{{\rangle}}
\def\B{\mathcal{B}}

\def\n{\nabla}

\def\a{\alpha}

\def\be{\begin{equation} }
\def\ee{\end{equation} }

\newcommand\blfootnote[1]{\begingroup
\renewcommand\thefootnote{}\footnote{#1}
\addtocounter{footnote}{-1}
\endgroup}
\begin{document}

\title{A topological sphere theorem for submanifolds 
of the hyperbolic space}
\author{M. Dajczer and Th. Vlachos}
\date{}
\maketitle

\begin{abstract}
We identify as topological spheres those complete submanifolds 
lying with any codimension in hyperbolic space whose Ricci 
curvature satisfies a lower bound contingent solely upon the 
length of the mean curvature vector of the immersion.
\end{abstract}

\blfootnote{\textup{2020} \textit{Mathematics Subject 
Classification}: 53C20, 53C40.}

There are numerous papers that characterize the topology 
of compact submanifolds in space forms of nonnegative 
sectional curvature under pinching assumptions that 
encompass both intrinsic and extrinsic data. The former 
are given in terms of some metric curvature, while the 
latter incorporates concepts derived from the second 
fundamental form of the submanifold, quite often emphasizing 
its norm. Most of these papers are \cite{DV1}, \cite{DV2}, 
\cite{OV}, \cite{V1}, \cite{V2}, \cite{X}, 
\cite{XG}, \cite{XHG}, \cite{XT} and \cite{XLG}.

In rather stark contrast, the situation diverges significantly 
when the ambient space form features negative sectional curvature, 
a scenario addressed by this paper. The papers we have been 
able to find pertaining  to this case, namely \cite{FX}, \cite{SX} 
and \cite{V3}, do not offer results related to the one presented 
here.\vspace{1ex}

Let $f\colon M^n\to\Hy^{n+m}_c$, $n\geq 4$, be an isometric 
immersion with codimension $m$ of a complete $n$-dimensional
Riemannian manifold into the hyperbolic space of constant sectional 
curvature $c<0$. Let $\Ric_M$ stand for the (not normalized) 
Ricci curvature of $M^n$ and denote the norm of the (normalized) 
mean curvature vector field $\mathcal{H}$ by $H$. 

\begin{theorem}\label{thm1}
Let $f\colon M^n\to\Hy_c^{n+m}$, $n\geq 4$, be an isometric 
immersion of a complete Riemannian manifold. If at any point 
it holds that
\be\label{1}
\Ric_M\geq (n-4)c+(n-2)H^2\tag{$\ast$}
\ee 
then $M^n$ is homeomorphic to a sphere $\Sf^n$.
\end{theorem}

 A well-known result due to Hamilton \cite{Ha} gives that 
for dimension $n=3$ the submanifold is diffeomorphic 
to a spherical space form. 

When $M^n$ possesses the topological structure of a sphere, 
a conjecture for the weaker bound $(n-2)(c+H^2)$ under the 
assumption $c+H^2\geq 0$ have been put forth by Xu and Gu \cite{XG}.
It proposes that the submanifold should not merely be topologically 
equivalent but diffeomorphic to a sphere. In our case this holds true 
for dimensions $n=5,6,12,56,61$ as in these cases it has been 
established that the differentiable structure is unique; see 
Corollary $1.15$ in \cite{WX}.\vspace{1ex}

There are plenty of compact submanifolds in the hyperbolic 
space that satisfy \eqref{1} strictly at any point. 
In fact, this is the case for any totally umbilical 
$n$-dimensional submanifold, being the inequality strict 
at any point only if $H>\sqrt{-3c}$. Notice that it will 
persist in its strict form after subjecting a totally 
umbilical submanifold to a sufficiently small smooth 
deformation.

Finally, notice that the theorem extends its applicability 
to compact submanifolds within both Euclidean space and the 
round sphere, as these are umbilical submanifolds of the 
hyperbolic space. However, in such cases the assumption 
regarding the Ricci curvature is more restrictive 
compared to that stipulated in \cite{DV2}.

\section{The pinching condition}

The following result for complete simply connected spaces forms 
of sectional curvature  $c>0$ has been proved by Lawson 
and Simons \cite{1}, and then by Xin \cite{X} for $c=0$ when 
strict inequality holds in \eqref{1} at any point. Elworthy 
and Rosenberg \cite[p.\ 71]{ER} observed that the result 
still holds by only requiring the bound to be strict at 
some point of the submanifold. In case $c<0$ the result was 
proved by Fu and Xu \cite{FX} if strict inequality holds 
at any point. But the observation by Elworthy and Rosenberg 
also applies to this case.

\begin{theorem}\label{ls} 
Let $f\colon M^n\to\Q_c^{n+m}$, $n\geq 4$
be an isometric immersion of a compact manifold and  $p$ an integer
such that $1\leq p\leq n-1$. Assume that at any point $x\in M^n$ 
and for \emph{any} orthonormal basis $\{e_j\}_{1\leq j\leq n}$ 
of $T_xM$ the second fundamental form 
$\alpha_f\colon TM\times TM\to N_fM$  satisfies
\be
\Theta_p\!=\!\sum_{i=1}^p\sum_{j=p+1}^n\!\!\big(2\|\a_f(e_i,e_j)\|^2
\!-\<\a_f(e_i,e_i),\a_f(e_j,e_j)\>\big)
\!\leq p(n-p\,{\rm{sign}}(c))c \tag{$\#$}.
\ee
If the inequality $(\#)$ is strict at some point $x\in M^n$ and 
for any orthonormal basis of $T_xM$, then there are no stable 
$p$-currents and the homology groups satisfy that
$H_p(M^n;\mathbb{Z})=H_{n-p}(M^n;\mathbb{Z})=0$.
\end{theorem}

In the sequel we will make use of the following lemma.

\begin{lemma}\label{fi}
Let $f\colon M^n\to\Hy_c^{n+m}$, $n\geq 3$, be an isometric 
immersion that  satisfies at $x\in M^n$ the inequality \eqref{1}. 
Then for the traceless part of the second fundamental form 
$\phi=\a_f-\<\,,\,\>\mathcal H$ at $x$ we have 
$\frac{1}{n}\|\phi\|^2\leq H^2+3c$.
Hence $H\geq\sqrt{-3c}$ and if $H=\sqrt{-3c}$ then $f$ is 
totally umbilical at $x$.
\end{lemma}

\proof From \eqref{1} the scalar curvature satisfies
$\tau\geq n(n-4)c+n(n-2)H^2$.
On the other hand, the Gauss equation gives that 
\be\label{tau}
\tau=n(n-1)c+n^2H^2-S,
\ee
where $S$ is the norm of the second fundamental form.  
Then $S\leq 3nc+2nH^2$. Since 
$\|\phi\|^2=S-nH^2$ we have the desired inequality.
\qed\vspace{2ex}

Recall that a vector in the normal space $\eta\in N_fM(x)$ 
at $x\in M^n$ is named a \emph{Dupin principal normal} of 
$f\colon M^n\to\Hy^{n+m}_c$ at $x\in M^n$
if the associated tangent vector subspace
$$
E_\eta(x)=\left\{X\in T_xM\colon\alpha_f(X,Y)
=\<X,Y\>\eta\;\,\text{for all}\;\,Y\in T_xM\right\}
$$
is at least two dimensional. The dimension of $E_\eta(x)$ 
is the \emph{multiplicity} of $\eta$.
\vspace{1ex}

The proof of the following result is inspired by computations
given by Xu and Gu in \cite{XG} and more recently by us in 
\cite{DV1}.

\begin{proposition}\label{prop}
Let $f\colon M^n\to\Hy_c^{n+m}$, $n\geq 4$, be an isometric 
immersion satisfying the inequality $(*)$ at $x\in M^n$.
Then at $x\in M^n$ the following assertions hold:
\vspace{1ex}

\noindent $(i)$ The inequality $(\#)$ is satisfied for any 
integer $2\leq p\leq n/2$ and for any orthonormal basis of $T_xM$. 
Moreover, if the inequality $(*)$ is strict or if $p<n/2$ then 
also $(\#)$ is strict for any orthonormal basis of $T_xM$. 
\vspace{1ex}

\noindent $(ii)$ Assume that equality holds in $(\#)$ for 
a certain integer $2\leq p\leq n/2$ and an orthonormal basis 
$\{e_j\}_{1\leq j\leq n}$ of $T_xM$. Then $n=2p$ and 
$$
\Ric_M(X)=(n-4)c+(n-2)H^2\;\; \text{for any unit}\;\; 
X\in T_xM.
$$
Moreover, we have:
\begin{itemize}
\item[(a)] If $n\geq 6$ then either $f$ is totally 
umbilical with $H=\sqrt{-3c}$ or we have that $H>\sqrt{-3c}$  
and there are distinct Dupin principal normals 
$\eta_1$ and $\eta_2$ such that 
$E_{\eta_1}=\spa\{e_1,\dots,e_p\}$ and 
$E_{\eta_2}=\spa\{e_{p+1},\dots,e_n\}$.
\item[(b)] If $n=4$ there are normal vectors $\eta_j$, 
$j=1,2$, such that 
\be\label{Vj}
\pi_{V_j}\circ A_\xi|_{V_j}=\<\xi,\eta_j\>I\;\,
\text{for any}\;\,\xi\in N_f(x)
\ee
where $V_1=\spa\{e_1,e_2\}$, $V_2=\spa\{e_3,e_e\}$ and 
$\pi_{V_j}\colon T_xM\to V_j$ are the projections.
\end{itemize}
\end{proposition}

\proof Recall that the Gauss equation of $f$ yields that 
the Ricci curvature  for any unit vector $X\in T_xM$ is
given by
\be\label{ric}
\Ric_M(X)=(n-1)c+\sum_{\a=1}^m(\mathrm{tr}A_\a)
\<A_\a X,X\>-\sum_{\a=1}^m\|A_{\a}X\|^2, 
\ee
where the $A_{\a}$, $1\leq\a\leq m$, stand for the shape 
operators of $f$ associated to an orthonormal basis 
$\{\xi_\alpha\}_{1\leq\alpha\leq m}$ of the normal vector 
space $N_fM(x)$.

From now on, the basis  
$\{\xi_\alpha\}_{1\leq\alpha\leq m}$ is taken such that
$\mathcal{H}(x)=H(x)\xi_1$ when $H(x)\neq 0$. For a given 
orthonormal basis $\{e_j\}_{1\leq j\leq n}$ of $T_xM$, we 
denote for simplicity $\a_{ij}=\alpha_f(e_i,e_j)$, 
$1\leq i,j\leq n$. Then, we have
\begin{align}\label{a}
\Theta_p&=\;2\sum_{i=1}^{p}\sum_{j=p+1}^n\|\a_{ij}\|^2-n
\sum_{i=1}^p\<\a_{ii},\mathcal{H}\>
+\sum_{i,j=1}^p\<\a_{ii},\a_{jj}\>\nonumber\\
=&\;2\sum_{i=1}^p\sum_{j=p+1}^n\sum_\a\<A_{\alpha }e_i,e_j\>^2-nH
\sum_{i=1}^p\<A_1e_i,e_i\>
+\sum_\a\big(\sum_{i=1}^p\<A_\alpha e_i,e_i\>\big)^2\nonumber\\
\leq&\; 2\sum_{i=1}^p\sum_{j=p+1}^n\sum_\a\<A_{\alpha }e_i,e_j\>^2-nH
\sum_{i=1}^p\<A_1e_i,e_i\>
+p\sum_\a\sum_{i=1}^p\<A_{\alpha}e_i,e_i\>^2,
\end{align}
where the inequality part was obtained using the Cauchy-Schwarz 
inequality
\be\label{a-}
\big(\sum_{i=1}^p\<A_\alpha e_i,e_i\> 
\big)^2\leq p\sum_{i=1}^p\<A_{\alpha}e_i,e_i\> ^2. 
\ee
Since $p\geq 2$ by assumption, then 
\begin{align}\label{b-}
2\sum_{i=1}^p&\sum_{j=p+1}^n\sum_\a\<A_{\alpha}e_i,e_j\>^2
+p\sum_{i=1}^p\sum_\a\<A_{\alpha}e_i,e_i\>^2\nonumber\\
&\leq p\sum_{i=1}^p\sum_{j=p+1}^n
\sum_\a\<A_{\alpha }e_i,e_j\>^{2}+p\sum_{i=1}^p\sum_\a
\<A_{\alpha}e_i,e_i\>^2\nonumber\\
&\leq p\sum_{i=1}^p\sum_\a\|A_{\alpha}e_i\|^2 
\end{align}
and thus \eqref{a} implies that
$$
\Theta_p\leq p\sum_{i=1}^p\sum_\a\|A_{\alpha}e_i\|^2
-nH\sum_{i=1}^p\<A_1e_i,e_i\>.
$$
Setting $\varphi=A_1-HI$ and using \eqref{ric}, we obtain 
\begin{align}\label{first}
\Theta_p
&\leq p\sum_{i=1}^p\left((n-1)c-\Ric_M(e_i)\right)
+(p-1)nH\sum_{i=1}^p\<A_1e_i,e_i\>\nonumber\\
&= p\sum_{i=1}^p\left((n-1)(c+H^2)-\Ric_M(e_i)\right)
-p(n-p)H^2\nonumber\\
&\;+(p-1)nH\sum_{i=1}^p\<\varphi e_i,e_i\>.
\end{align}
Then
\be\label{d}
\Theta_p
\leq p^2\big((n-1)(c+H^2)-\Ric_M^{\text{min}}(x)\big)-p(n-p)H^2
+(p-1)nH\sum_{i=1}^p\<\varphi e_i,e_i\>
\ee
where 
$$
\Ric_M^{\text{min}}(x)
=\min\left\{\Ric_M(X)\colon X\in T_xM,\|X\|=1\right\}.
$$
\vspace{1ex}

We have that 
\be\label{min}
(n-1)(c+H^2)\geq\Ric_M^{\text{min}}(x)
\ee
and that equality holds if $f$ is totally umbilical at $x\in M^n$.
In fact, it follows from \eqref{tau} that
\be\label{sca}
S\leq n(n-1)c+n^2H^2-n\Ric_M^{\text{min}}(x).
\ee
Therefore, 
\begin{align*}
(n-1)(c+H^2)-&\Ric_M^{\text{min}}(x)\geq\frac{1}{n}(S-nH^2)
=\frac{1}{n}\|\phi\|^2,
\end{align*}
and \eqref{min} follows.

From \eqref{min} and having that $p\leq n/2$, then
\be\label{re}
p^2\big((n-1)(c+H^2)-\Ric_M^{\text{min}}(x)\big)\leq p(n-p)
\big((n-1)(c+H^2)-\Ric_M^{\text{min}}(x)\big).
\ee
Therefore, from \eqref{d} it follows the estimate
\be\label{c1}
\Theta_p
\leq p(n-p)\big((n-1)(c+H^2)-\Ric_M^{\text{min}}(x)-H^2\big)
+(p-1)nH\sum_{i=1}^p\<\varphi e_i,e_i\>.
\ee

Next, we obtain a second estimate of 
$$
\Theta_p=\;\sum_\a\Big(2\sum_{i=1}^p\sum_{j=p+1}^n\<A_\a e_i,e_j\>^2
-\sum_{i=1}^p\<A_\a e_i,e_i\>\sum_{j=p+1}^n\<A_\a e_j,e_j\>\Big).
$$
Decomposing
\begin{align*}
&\sum_{i=1}^p\<A_\a e_i,e_i\>\sum_{j=p+1}^n\<A_\a e_j,e_j\>\\
=&\frac{n-p}{n}\sum_{i=1}^p\<A_\a e_i,e_i\>
\sum_{j=p+1}^n\<A_\a e_j,e_j\>
+  \frac{p}{n}\sum_{i=1}^p\<A_\a e_i,e_i\>
\sum_{j=p+1}^n\<A_\a e_j,e_j\>,
\end{align*}
we have
\begin{align*}
&\Theta_p
=\sum_\a\Big(2\sum_{i=1}^p\sum_{j=p+1}^n\<A_\a e_i,e_j\>^2
-\frac{n-p}{n}\trace A_\a\sum_{i=1}^p\<A_\a e_i,e_i\>\\
&+\frac{n-p}{n}\big(\sum_{i=1}^p\<A_\a e_i,e_i\>\big)^2
-\frac{p}{n}\trace A_\a\sum_{j=p+1}^n\<A_\a e_j,e_j\>
+\frac{p}{n}\big(\sum_{j=p+1}^n\<A_\a e_j,e_j\>\big)^2\Big).
\end{align*}
Since $p(n-p)/n\geq 1$, we obtain using the Cauchy-Schwarz 
inequality that
\begin{align*}\label{a1}
\Theta_p
\leq&\sum_\a\Big(2\sum_{i=1}^p\sum_{j=p+1}^n\<A_\a e_i,e_j\>^2
-\frac{n-p}{n}\trace A_\a\sum_{i=1}^p\<A_\a e_i,e_i\>\\
&+\;\frac{p(n-p)}{n}\sum_{i=1}^p\<A_\a e_i,e_i\>^2
-\frac{p}{n}\trace A_\a\sum_{j=p+1}^n\<A_\a e_j,e_j\>\\
&+\;\frac{p(n-p)}{n}\sum_{j=p+1}^n\<A_\a e_j,e_j\>^2\Big)\\
\leq&\;\frac{p(n-p)}{n}S-(n-p)H\sum_{i=1}^p\<A_1e_i,e_i\>
-pH\sum_{j=p+1}^n\<A_1e_j,e_j\>.
\end{align*}
It follows that 
$$
\Theta_p
\leq\frac{p(n-p)}{n}S-2p(n-p)H^2-(n-2p)H
\sum_{i=1}^p\<\varphi e_i,e_i\>.
$$
Then we have from \eqref{sca} that
\be\label{a2}
\Theta_p
\leq p(n-p)\big((n-1)(c+H^2)-\Ric_M^{\text{min}}(x)-H^2\big)
-(n-2p)H\sum_{i=1}^p\<\varphi e_i,e_i\>.
\ee
By computing $(n-2p)\times$\eqref{c1}$+n(p-1)\times$\eqref{a2}, 
we obtain
\be\label{ac1}
\Theta_p\leq p(n-p)\left((n-1)c+(n-2)H^2-\Ric_M^{\text{min}}(x)\right).
\ee
It follows from \eqref{ac1} using $(*)$ that
\be\label{ac2}
\Theta_p-p(n+p)c\leq 2p(n-2p)c,
\ee
and the inequality $(\#)$ has been proved. Clearly, if the 
inequality $(*)$ is strict or if $p<n/2$ then $(\#)$ is strict, 
and this completes the proof of part $(i)$.
\vspace{1ex}

We prove part $(ii)$. From part $(i)$ we have that 
$n=2p$ and that all the inequalities that range from 
\eqref{a} to \eqref{d} as well as the ones from \eqref{re} 
to \eqref{ac2} become equalities.

From \eqref{a-} we obtain 
\be\label{r}
\<A_\a e_i,e_i\>=\rho_\a\;\,\text{for all}
\;1\leq i\leq p,\;\,1\leq\a\leq m.
\ee
We have that  \eqref{b-} gives
\be\label{k}
(p-2)\<A_\a e_i,e_j\>=0\;\,\text{for all}
\;\,1\leq i\leq p,\;p+1\leq j\leq n,\;1\leq\a\leq m, 
\ee
and
\be\label{ii}
\<A_\a e_i,e_{i'}\>=0\;\,\text{for all}\;\,1\leq i\neq i' 
\leq p,\;1\leq\a\leq m.
\ee
From \eqref{first} and \eqref{d} we have 
$\Ric_M(e_i)=\Ric_M^{\text{min}}(x)$. Then 
\eqref{ac1} and \eqref{ac2} give 
\be\label{ricb}
\Ric_M(e_i)=\Ric_M^{\text{min}}(x)=(n-4)c+(n-2)H^2
\;\,\text{for all}\;\,1\leq i\leq p.
\ee
Since $n=2p$ then equality also holds in $(\#)$ for
the reordered orthonormal basis 
$\{e_{p+1},\dots,e_n, e_1,\dots,e_p\}$ of $T_xM$. 
Therefore, we also have that
\be\label{r1}
\<A_\a e_j,e_j\>=\mu_\a\;\,\text{for all}
\;\,p+1\leq j\leq n,\;\,1\leq\a\leq m,
\ee
\be\label{jj}
\<A_\a e_j,e_{j'}\>=0\;\,\text{for all}\;\,p+1\leq j\neq j' 
\leq n,\;1\leq\a\leq m,
\ee
and 
\be\label{ricbj}
\Ric_M(e_j)=\Ric_M^{\text{min}}(x)=(n-4)c+(n-2)H^2
\;\,\text{for all}\;\,p+1\leq j\leq n.
\ee
Hence, we obtain from \eqref{ricb} and \eqref{ricbj}  
that 
\be\label{ein}
\Ric_M(X)=(n-4)c+(n-2)H^2\;\,\text{for any unit}\;\,X\in T_xM.
\ee

In particular, it follows at $x$ from \eqref{r}, \eqref{k}, 
\eqref{ii}, \eqref{r1} and \eqref{jj} that the vectors
$\eta_1=\sum_\a\rho_\a\xi_\a$ and 
$\eta_2=\sum_\a\mu_\a\xi_\a$ are Dupin principal normals 
with 
$
E_{\eta_1}=\spa\{e_1,\dots,e_p\}\;\;
{\text {and}}\;\;E_{\eta_2}=\spa\{e_{p+1},\dots,e_n\}.
$
If $\eta_1=\eta_2$, then $f$ at $x$ is totally umbilical and 
equality holds in \eqref{min}. This combined with \eqref{ein} 
yields $H=\sqrt{-3c}$ at $x$. 
If otherwise, then Lemma \ref{fi} gives at $x$ that $H>\sqrt{-3c}$ 
and that $\eta_1$ and $\eta_2$ are distinct Dupin principal 
normals, and this concludes the proof of part~$(a)$. 

Finally, if $n=4$ then for any $\xi\in N_f(x)$ we have 
\eqref{Vj} where $\eta_1=\sum_\a\rho_\a\xi_\a$ and 
$\eta_2=\sum_\a\mu_\a\xi_\a$, and part $(ii)$ has also 
been proved.\qed

\section{The proof of Theorem \ref{thm1}}

For the proof of Theorem~\ref{thm1} we initially establish 
a topological result.

\begin{lemma}\label{p1}
Let $f\colon M^n\to\Hy_c^{n+m}$, $n\geq 4$, 
be an isometric immersion of a compact manifold satisfying
\be\label{p=1}
\Ric_M> \frac{1}{n+2}\big((n^2-n-4)c+n(n-1)H^2\big).
\ee
Then $\pi_1(M^n)=0$ and $H_{n-1}(M^n,\Z)=0$.
\end{lemma}

\proof From \eqref{tau} and \eqref{p=1} it follows that 
\be\label{phi}
\|\phi\|^2\leq \frac{2n}{n+2}\big((n+1)c+(n-1)H^2\big).
\ee

Let $\{e_i\}_{1\leq i\leq n}$ and $\{\xi_\a\}_{1\leq\a\leq m}$
be orthonormal tangent and normal bases  at $x\in M^n$.
Using \eqref{ric} we obtain that
\begin{align*}
\sum_{j=2}^n\big(2&\|\a_{1j}\|^2
-\<\a_{11},\a_{jj}\>\big)\\
=&\;2\sum_\a\sum_{j=2}^n\<A_\a e_1,e_j\>^2
-\sum_\a\<A_\a e_1,e_1\>\sum_{j=2}^n\<A_\a e_j,e_j\>\\
=&\;\sum_\a\sum_{j=2}^n\<A_\a e_1,e_j\>^2
-\sum_\a\trace A_\a\<A_\a e_1,e_1\>
+\sum_\a\|A_\a e_1\|^2\\
=&\;\sum_{j=2}^n\|\phi( e_1,e_j)\|^2+(n-1)c-\Ric_M(e_1).
\end{align*}
Then this together with \eqref{p=1} and \eqref{phi} give
$$
\sum_{j=2}^n\big(2\|\a_{1j}\|^2-\<\a_{11},\a_{jj}\>\big)
\leq\;\frac{1}{2}\|\phi\|^2+(n-1)c-\Ric_M(e_1)<(n+1)c.
$$
Hence, by Theorem \ref{ls} there are no stable $1$-currents 
on $M^n$ and therefore $H_1(M^n,\Z)=H_{n-1}(M^n,\Z)=0$. 
Since in each nontrivial free homotopy class there is a 
length minimizing curve, we conclude that $\pi_1(M^n)=0$.
\vspace{2ex}\qed

\noindent\emph{Proof of Theorem \ref{thm1}:} 
We have from Lemma \ref{fi}  that $H\geq\sqrt{-3c}$.  
It then follows from \eqref{1} that $\Ric_M\geq-2(n-1)c$ 
and hence  $M^n$ is compact by the classical Bonnet-Myers 
theorem. Moreover, since
$$
\Ric_M\geq (n-4)c+(n-2)H^2
>\frac{1}{n+2}\big((n^2-n-4)c+n(n-1)H^2\big)
$$
then Lemma \ref{p1} yields that $M^n$ is simply connected 
and  $H_{n-1}(M^n,\Z)=0$.

According to part $(i)$ of Proposition \ref{prop} the 
inequality $(\#)$ is satisfied at any point of $M^n$ 
for any $2\leq p\leq n/2$ and for any orthonormal tangent 
basis at that point. 
We argue that the homology groups satisfy
\be\label{Hom}
H_p(M^n;\Z)=0=H_{n-p}(M^n;\Z)\;\,{\text {for all}}
\;\,2\leq p\leq n/2. 
\ee
Suppose to the contrary that \eqref{Hom} does 
not hold. Consider the nonempty set
$$
P=\left\{2\leq p\leq n/2: H_p(M^n;\Z)
\neq0\;\,{\text {or}}\;\,H_{n-p}(M^n;\Z)\neq0\right\}
$$
and denote $k=\max P$. Hence $H_k(M^n;\Z)\neq 0$ or 
$H_{n-k}(M^n;\Z)\neq 0$. By Theorem \ref{ls} at any point 
$x\in M^n$ there is an orthonormal tangent basis
such that equality holds in $(\#)$ for $p=k$. Moreover, 
we have from part $(ii)$ of 
Proposition \ref{prop} that \eqref{ein} holds. In this 
situation it is well-known that $M^n$ is an Einstein 
manifold and, in particular, it follows that $H$ is a 
positive constant.\vspace{1ex} 

We need to differentiate between two cases based on 
the dimension of the submanifold.
\vspace{1ex}

\noindent\emph{Case $n\geq 6$}.
Part $(ii)$ of Proposition \ref{prop} yields $k= n/2$. 
We argue that $H>\sqrt{-3c}$. If we have otherwise, then 
the submanifold is totally umbilical by part $(ii)$ of 
Proposition \ref{prop}. Hence, we have from the Gauss 
equation that $M^n$ has constant sectional curvature $-2c$. 
But then $M^n$ would be isometric to a round sphere, and 
this contradicts our assumption that $H_k(M^n;\Z)\neq 0$. 

Since $H>\sqrt{-3c}$, according to 
part $(ii)$ of Proposition \ref{prop} there are smooth 
Dupin principal normal vector fields $\eta_1$ and $\eta_2$ of 
multiplicity $k$ and corresponding smooth distributions 
$E_1$ and $E_2$. Let $\{X_\ell\}_{1\leq\ell\leq n}$ be a 
smooth local orthonormal frame satisfying that
$E_1=\spa\left\{X_1,\dots,X_k\right\}$ and 
$E_2=\spa\left\{X_{k+1},\dots,X_n\right\}$. Then 
$\a_f(X_i,X_i)=\eta_1$ if $1\leq i\leq k$ and 
$\a_f(X_j,X_j)=\eta_2$ if $k+1\leq j\leq n$.
\vspace{1ex}

If follows from the Gauss equation that 
$$
\Ric_M(X)=(n-1)c\|X\|^2+n\<\mathcal{H},\a_f(X,X)\>-I\!I\!I(X)
\;\,\text{for any}\;\,X\in\mathcal X(M),
$$
where $I\!I\!I(X)=\sum_{\ell=1}^n\|\a_f(X,X_\ell)\|^2$ is 
the so called third fundamental form of $f$. Since 
$\mathcal H=(\eta_1+\eta_2)/2$, then
\be\label{H}
4H^2=\|\eta_1\|^2+\|\eta_2\|^2+2\<\eta_1,\eta_2\>.
\ee
Moreover, we have for $1\leq i\leq k$ that
$$
I\!I\!I(X_i)=\sum_{\ell=1}^n\|\a(X_\ell,X_i)\|^2
=\|\eta_1\|^2
$$
and
\begin{align*}
(n-4)c+(n-2)H^2&=\Ric_M(X_i)
=(n-1)c+n\<\mathcal H,\a(X_i,X_i)\>-I\!I\!I(X_i)\\
&=(n-1)c+k\<\eta_1+\eta_2,\eta_1\>-\|\eta_1\|^2.
\end{align*}
Thus
\be\label{ric1}
(n-4)c+(n-2)H^2=(n-1)c+(k-1)\|\eta_1\|^2+k\<\eta_1,\eta_2\>.
\ee
Arguing similarly for $k+1\leq j\leq n$, we obtain
\be\label{ric2}
(n-4)c+(n-2)H^2=(n-1)c+(k-1)\|\eta_2\|^2+k\<\eta_1,\eta_2\>.
\ee
It follows from \eqref{ric1} and \eqref{ric2} that
$\|\eta_1\|=\|\eta_2\|$, and hence \eqref{H} 
becomes
\be\label{H1}
2H^2=\|\eta_1\|^2+\<\eta_1,\eta_2\>.
\ee
Combining \eqref{ric1} with \eqref{H1} gives
\be\label{12}
\<\eta_1,\eta_2\>=-3c. 
\ee
Then, we conclude from \eqref{H1} that 
\be\label{1122}
\|\eta_1\|^2=\|\eta_2\|^2=2H^2+3c>0.
\ee

The Codazzi equation for $f$ is easily seen to yield
\be\label{cod}
\<\nabla_XY,Z\>(\eta_i-\eta_j)=\<X,Y\>\nabla_Z^\perp\eta_i
\;\,\text{if}\;\,i\neq j
\ee
for any $X,Y\in E_i,Z\in E_j$.
Using \eqref{12} and \eqref{1122} then \eqref{cod} gives
$$
2\<\nabla_XY,Z\>(H^2+3c)
=\<X,Y\>\<\nabla_Z^\perp\eta_i,\eta_i\>=0
$$
for all $X,Y\in E_i$ and $Z\in E_j,\;i\neq j$, that is, 
the distributions $E_1$ and $E_2$ are totally geodesic. 
Being simply connected, then de Rham theorem gives that 
$M^n$ is a Riemannian product $M_1^k\times M_2^k$ such 
that $TM_j^k=E_j$, $j=1,2$. It follows from the Gauss 
equation that the manifolds $M_1^k$ and $M_2^k$ have both 
constant sectional curvature $2H^2+4c$. But then the Ricci 
curvature of $M^n=M_1^k\times M_2^k$ is $(n-2)(H^2+2c)$, 
which is in contradiction with \eqref{ein}.
\vspace{1ex}

\noindent\emph{Case $n=4$}. We have that $k=2$ and 
$H_2(M^4;\Z)\neq 0$. Since $\Ric_M=2H^2$ then $\tau=8H^2$ 
and hence \eqref{tau} gives $S=12c+8H^2$ or, equivalently, that 
$\|\phi\|^2=12c+4H^2$. It then follows from Proposition 
$16$ in \cite{OV} that the Bochner operator 
$\B^{[2]}\colon\Omega^2(M^4)\to\Omega^2(M^4)$, 
a certain symmetric endomorphism of the bundle of $2$-forms 
$\Omega^2(M^4)$, satisfies for any $\omega\in\Omega^2(M^4)$ 
the inequality
\be\label{ineq}
\<\B^{[2]}\omega,\omega\>
\geq \big((4(c+H^2)-\|\phi\|^2\big)\|\omega\|^2
=-8c\|\omega\|^2.\nonumber
\ee
Hence $\B^{[2]}$ is positive definite. 
 
We claim that the second Betti number $\beta_2(M^4)$ of 
the manifold vanishes. If otherwise, then there would exist 
a nonzero harmonic $2$-form $\omega\in\Omega^2(M^4)$. By the 
Bochner-Weitzenb\"ock formula the Laplacian of 
$\omega$ is given by
$$
0=\Delta\omega=\n^*\n\omega+\B^{[2]}\omega,
$$
where $\n^*\n$ is the rough Laplacian. From this we obtain 
\be\label{boch.formula}
\|\n\omega\|^2+\<\B^{[2]}\omega,\omega\>
+\frac{1}{2}\,\Delta \|\omega\|^2=0.\nonumber
\ee 
Then the maximum principle and the fact that $\B^{[2]}$ is 
positive definite imply that $\omega=0$, which proves the
claim. 

From the claim, we have that $H_2(M^4;\Z)$ is a nontrivial 
torsion group and  Poincar\'e duality gives that the torsion 
of $H^2(M^4;\Z)$ is isomorphic to $H_2(M^4;\Z)$. On the other 
hand, the universal coefficient theorem of cohomology yields 
that the torsion subgroups of $H^2(M^4;\Z)$
and $H_1(M^4;\Z)$ are isomorphic 
(cf.\ \cite[p.\ 244 Corollary 4]{Sp}). Since $M^4$ is simply 
connected, we have that $H_1(M^n;\Z)=0$ and thus $H^2(M^4;\Z)$ 
is torsion free. This is a contradiction and completes 
the proof that \eqref{Hom} holds.

Hence $M^n$ is a simply connected homology sphere and it follows 
from the Hurewicz isomorphism theorem that $M^n$ is a homotopy 
sphere. Finally, the resolution of the generalized Poincar\'{e} 
conjecture gives that $M^n$ is homeomorphic to $\Sf^n$. 
\vspace{2ex}\qed

\noindent Marcos Dajczer is  partially supported 
by the grant PID2021-124157NB-I00 funded by 
MCIN/AEI/10.13039/501100011033/ `ERDF A way of making Europe',
Spain, and are also supported by Comunidad Aut\'{o}noma de la 
Regi\'{o}n de Murcia, Spain, within the framework of the Regional 
Programme in Promotion of the Scientific and Technical Research 
(Action Plan 2022), by Fundaci\'{o}n S\'{e}neca, Regional Agency 
of Science and Technology, REF, 21899/PI/22.
\bigskip

Theodoros Vlachos thanks the Department of Mathematics 
of the University of Murcia where part of this work was 
done for its cordial hospitality during his visit. 
He was supported by the grant PID2021-124157NB-I00 
funded by MCIN/AEI/10.13039/501100011033/ `ERDF 
A way of making Europe', Spain.

\noindent Marcos Dajczer\\
Departamento de Matemáticas\\ 
Universidad de Murcia, Campus de Espinardo\\ 
E-30100 Espinardo, Murcia, Spain\\
e-mail: marcos@impa.br
\bigskip

\noindent Theodoros Vlachos\\
Department of Mathematics\\
University of Ioannina\\
45110 Ioannina, Greece\\
e-mail: tvlachos@uoi.gr

\begin{thebibliography}{lll}

\bibitem{DT} Dajczer, M. and Tojeiro, R., 
``Submanifold theory beyond an introduction". 
Springer, 2019.

\bibitem{DV1} Dajczer, M. and Vlachos, T.,
\textit{Ricci pinched compact submanifolds in spheres},
Preprint.

\bibitem{DV2} Dajczer, M. and Vlachos, T.,
\textit{Ricci pinched compact submanifolds in space forms}. 
Preprint.

\bibitem{ER} Elworthy, K. and Rosenberg, S., 
\textit{Homotopy and homology vanishing theorems and 
the stability of stochastic flows}, 
Geom. Funct. Anal. \textbf{6} (1996), 51--78.

\bibitem{FX} Fu, H. and Xu, H.,
\textit{Vanishing and topological sphere theorems for 
submanifolds in a hyperbolic space}, 
Internat. J. Math. \textbf{19} (2008), 811--822.

\bibitem{Ha} Hamilton, R., 
\textit{Three-manifolds with positive Ricci curvature}
J. Differential Geom. \textbf{17} (1982), 255--306.

\bibitem{1} Lawson, B. and Simons, J.,
\textit{On stable currents and their application to global 
problems in real and complex geometry},  
Ann. of Math. \textbf{98} (1973), 427--450.

\bibitem{OV} Onti, C. R. and Vlachos, Th.,
\textit{Homology vanishing theorems for pinched submanifolds}, 
J. Geom. Anal. \textbf{222} (2022), Paper 294, 34 pp.

\bibitem{SX} Shiohama, K. and Xu, H., 
\textit{The topological sphere theoremfor complete submanifolds}, 
Compositio Math. \textbf{107} (1997), 221--232.

\bibitem{Sp} Spanier, E., ``Algebraic topology",
McGraw-Hill Book Co., New York-Toronto, Ont.-London, 1966.

\bibitem{V1} Vlachos, Th.,
\textit{A sphere theorem  for odd-dimensional submanifolds of spheres},
Proc. Amer. Math. Soc. \textbf{130} (2002), 167--173.

\bibitem{V2} Vlachos, Th.,
\textit{Homology vanishing theorems for submanifolds},
Proc. Amer. Math. Soc. \textbf{135} (2007), 2607--2617.

\bibitem{V3} Vlachos, Th.,
\textit{Geometric and topological rigidity of pinched submanifolds},
Preprint.

\bibitem{WX} Wang, G. and Xu, Z.,
\textit{The triviality of the 61-stem in the stable
homotopy groups of spheres},
Annals of Math. \textbf{186} (2017), 501--580.

\bibitem{X} Xin, Y., 
\textit{An application of integral currents to the
vanishing theorems}, 
Sci. Sinica Ser. A \textbf{27} (1984), 233--241.

\bibitem{XG} Xu, H. and Gu, J.,
\textit{Geometric, topological and differentiable rigidity 
of submanifolds in space forms}, 
Geom. Funct. Anal. \textbf{23} (2013), 1684--1703.

\bibitem{XHG} Xu, H, Huang, F. and Zhao, E.,
\textit{Entao Differentiable pinching theorems for 
submanifolds via Ricci flow}, 
Tohoku Math. J. \textbf{67} (2015), 531--540.

\bibitem{XT} Xu, H. and Tian, L., 
\textit{A differentiable sphere theorem inspired by 
rigidity of minimal submanifolds}, 
Pacific J. Math. \textbf{254} (2011), 499--510.

\bibitem{XLG} Xu, H., Leng, Y. and Gu, J., 
\textit{Geometric and topological rigidity for compact 
submanifolds of odd dimension}, 
Sci. China Math. \textbf{57} (2014), 1525--1538. 

\end{thebibliography}
\end{document}